\documentclass[12pt]{{article}}
\usepackage{amsmath,amsfonts,amssymb,amsthm,graphics,amscd}
\begin{document}

\title{Remarks and improvement regarding the theorem of continuous dependence of solutions of F.D.E., using some new results for continuous convergences}
\author{E. Athanasiadou, C. Papachristodoulos\\
Department of Mathematics\\
University of Athens,\\
GR-15784, Panepistimiopolis, Athens, Greece\\
email: eathan@math.uoa.gr, christos.papachristodoulos@math.uoa.gr }
\date{}
\maketitle

\begin{abstract}
\noindent We improve the theorem on continuous dependence of solutions of functional differential equations (see J. Hale, Functional differential equations, theorem 5.1), using some new results on continuous convergences. Namely, we prove this theorem without the assumption of uniform boundeness for the sequence $(f_k)$. Also we restore some inconsistencies regarding the proof of this theorem.
\end{abstract}
\smallskip
\noindent \textbf{Keywords:} continuous convergence or $a$-convergence, exhaustive sequence, weak exhaustive sequence, operator corresponding to F.D.E. of retarded type, Caratheodory condition\\
\noindent \textbf{Mathematics Subject Classification:} ............

\section{Introduction}
Let $r\geq 0, \sigma \in \mathbb{R}$ and $C=C([-r,0],R^N)$, $C_{\sigma , r}=C([\sigma -r,\sigma +a],R^N)$ $(a>0)$ be the Banach spaces of continuous functions with domains $[-r,0]$, $[\sigma -r,\sigma +a]$ and values into $\mathbb{R}^N$, endowed with the supremum norm. For each $x\in C_{\sigma,r}$ and $t\in[\sigma,\sigma+a]$ we set $x_t\in C$ to be the function
$$
x_t(\theta)=x(t+\theta), -r\leq\theta\leq0.
$$
If $D$ is an open subset at $\mathbb{R}\times C$, $f: D\rightarrow \mathbb{R}^N$ is a given function and $\varphi\in C$ is a fixed function, we consider the following most general problem
\[\left.
  \begin{array}{lr}
   \dot{x}(t)=f(t,x_t), t\in(\sigma,\sigma +a)  \\
    x_{\sigma}=\varphi
  \end{array}
\right\} (P),
\]
which is known as functional differential equation of retarded type with initial value $\varphi$ and starting point $\sigma$. For the sequel, we assume that $f$  satisfies Caratheodory condition on $D$, that is,
\begin{itemize}
  \item[1st] $f(t,\varphi)$ is measurable in $t$ for each fixed $\varphi$
  \item[2nd] $f(t,\varphi)$ is continuous in $\varphi$ for each fixed $t$
  \item[3rd] For each $(t,\varphi)\in D$ there is a neighborhood (briefly nbd for the sequel) $V$ of $(t,\varphi)$ and a Lebesque intergable function $m$ such that
      $$
      |f(s,\psi)|\leq m(s), \mbox{ for all } (s,\psi)\in V
      $$
\end{itemize}

We remind that problem $(P)$ has a solution whenever $f$ satisfies Caratheodory condition and that $x$ is a solution of $(P)$, if there is $a>0$ such that $x\in C_{\sigma,r}, x_{\sigma}=\varphi$ and $\dot{x}(t)=f(t,x_t)$ almost every where on $[\sigma,\sigma+a)$.

The theorems on continuous dependence and generally on limiting behavior of solutions of the above problem (P) are heavily based on a kind of convergence named continuous convergence or $a$-convergence. This mode of convergence is also very important in Harmonic Analysis, Function spaces theory etc (Zygmund called this convergence as uniform convergence at a point, see [Z], page 58). In the recent years, some new results regarding continuous convergence have came in order (We present these results in section 2). Using these results, in section 3 we obtain improvement of the theorem on continuous dependence of solutions of the problem (P) ([H], theorem 5.1). Precisely, we prove this theorem without the strong requirement that the sequence $(f_k)$ is uniformly bounded and we give exact meaning on the convergence of the sequence of solutions with the aid of the generalized continuous convergence (definition 2.6). Finally we restore some inconsistencies regarding the proof of this theorem (Remarks 3.6)

\section{Some new results regarding continuous convergence}
Let $(X,d)(Y,p))$ be arbitrary metric spaces (In particular $X\subseteq R\times C$, $Y=R^N$ in our case) and $f_n, f:X\rightarrow Y$, $n=1,2,...$. We recall the following definitons (See [APP], [GP]).
\\\\
\textbf{Definitions 2.1:}\\
(a) We say that the sequence $(f_n)$ converges continuously to $f$ ($f_n\xrightarrow{{\mbox{cont}}} f$) iff, for each $x\in X$ and for each sequence $(x_n)$ in $X$ with $x_n\rightarrow x$ it holds that $f_n(x_n)\rightarrow f(x)$.
\\
(b) We say that the sequence $(f_n)$ is exhaustive, iff
\begin{eqnarray*}
\forall\;x\in X\;\forall \varepsilon>0\;\exists\;\delta=\delta(x,\varepsilon)>0\;\exists n_0=n_0(x,\varepsilon):\\ d(x,t)<\delta\Rightarrow\rho(f_n(x),f_n(t))<\varepsilon ,\mbox{ for } n\geq n_0.
\end{eqnarray*}
(c) We say that $(f_n)$ is weakly-exhaustive, iff,
$$
\forall\;x\in X\;\forall \varepsilon>0\;\exists \delta=\delta(x,\varepsilon)>0 : d(t,x)<\delta\Rightarrow\exists\; n_t\in N:\rho(f_n(x),f_n(t))<\varepsilon\mbox{ for } n\geq n_t
$$
Obviously if $(f_n)$ is exhaustive then $(f_n)$ is weakly-exhaustive. It is not hard to see that the inverse implications fails ([GP]).\\
Now, we formulate some new results on continuous convergence. The proofs of propositions 2.2, 2.3 below can be found in [GP] and the proofs of propositions 2.4, 2.5 in [APP].
\\\\
\textbf{Proposition 2.2:}\\
The following are equivalent\\
(a)$f_n \xrightarrow{\mbox{cont}} f$ \\
(b) $(f_n)$ converges pointwise to $f$ and $(f_n)$ is exhaustive.
\\\\
\textbf{Proposition 2.3:}\\
Suppose that $(f_n)$ converges pointwise to $f$. Then the following are equivalent.\\
(a)$f$ is continuous \\
(b) $(f_n)$ is weakly exhaustive.\\

We note that the functions $f_n$, $n=1,2,...$ need not to be continuous in the above theorem. Also as a corollary from proposition 2.2,2.3 we get  that continuous limit of any sequence of functions is necessarily continuous. With the next theorems we see how continuous convergence and uniform convergence are related.
\\\\
\textbf{Proposition 2.4:}\\
Suppose that
$$
f_n\xrightarrow{\mbox{cont}} f.
$$
Then $(f_n)$ converges uniformly to $f$ on compact subsets of $X$.
\\\par
In the inverse direction, the continuity of the functions $f_n$, $n=1,2,...$ is necessary (For details and concrete examples regarding the difference between continuous convergence and uniform local convergence see [APP]).
\\\\
\textbf{Proposition 2.5:}\\
Suppose that $\{f_n\}\subseteq C(X,Y):=\{f:X\rightarrow Y|f \mbox{ is continuous}\}$. If for each $x\in X$ there is a neiborhood $A$ of $x$ such that $(f_n)$ converges uniformly to $f$ on A, then
$$
f_n\xrightarrow{\mbox{cont}} f.
$$
In case that $X$ is locally compact and $\{f_n\}\subseteq C(X,Y)$ as a corollary from proposition 2.4, 2.5 we get that
$$
f_n\xrightarrow{\mbox{cont}} f\Leftrightarrow(f_n)\mbox{ converges uniformly on compacta to } f.
$$
Continuous convergence can be generalized for sequences of functions which have not the same domain (see [K]).
\\\\
\textbf{Proposition 2.6:}\\
Let $D_n \subseteq X, f_n: D_n\rightarrow Y, n=0,1,2,...$. We say that the sequence $(f_n)_{n=1,2,...}$ converges continuously to $f_0$ $(f_n\xrightarrow{\mbox{cont}}f_0)$, iff, $x_n\rightarrow x_0,\; n\rightarrow\infty,\;\; x_n\in D_n,\;\; n=0,1,2,...\Rightarrow f_n(x_n)\rightarrow f_0(x_0),\; n\rightarrow \infty$.
\\\\
\textbf{Remark 2.7:}\\
If $f_n\xrightarrow{\mbox{cont}} f_0$ according to the above generalized notion, then, by proposition 2.4, $(f_n)$ converges uniformly to $f_0$ on each compact subset $K$ of $X$ such that $K\subseteq D_n$, $n=0,1,2,...$.

\section{Continuous dependence}
In this paragraph we keep the notations introduced in $\S 1$ and we use the same terminology as in [H]. More precisely:
\begin{itemize}
  \item $\mathbb{N}$ is the set of positive integers.
  \item Let $\sigma\in\mathbb{R}, a>0$. For each $\varphi\in C$, by $\widetilde{\varphi}\in C_{\sigma,r}$ we denote the following function
      $$
      \widetilde{\varphi}_{\sigma}=\varphi \;\&\; \widetilde{\varphi}(\sigma+t)=\varphi(0) \mbox{ for } t\in[0,a].
      $$
  \item If $a,\beta>0, A(a,\beta)$ is the set of functions $\eta\in C_{0,r}=C([-r,a],\mathbb{R}^N)$ such that\\ $\eta_0=0\; \&\; ||\eta||<\beta$ ($||\;||$ denotes the supremum norm). Apparently, $A(a,\beta)$ is a closed, bounded, convex subset of $C_{0,r}$
  \item Suppose $\sigma\in\mathbb{R},\varphi\in C$ are fixed. Then we define the following operator corresponding to problem (P).
      $$
      T:A(a,\beta)\rightarrow C_{0,r}
      $$
      $$
      \eta\rightarrow T\eta:T\eta(t)=\left\{
                \begin{array}{ll}
                  \int_0^tf(\sigma+s,\widetilde{\varphi}_{\sigma+s},\eta_s)ds, t\in[0,a]\\
                  0, t\in[-r,0]
                \end{array}
              \right.
      $$
      where the function $f:D\subseteq \mathbb{R}\times C\rightarrow\mathbb{R}^N$ satisfies Caratheodory condition on $D$.
  \item The expression $\varphi_n\xrightarrow{||\;||}\varphi$, $n\rightarrow\infty$ means that the sequence of functions $\varphi_n$ converges uniformly to the function $\varphi$.
\end{itemize}
We point out for the sequel the following known facts
\\\\
\textbf{Fact 3.1:}\\
The following are equivalent.
\begin{itemize}
  \item[(a)] $x\in C_{\sigma ,r}=C([\sigma -r,\sigma+a],\mathbb{R}^N) \& x_{\sigma}=\varphi \in C$
  \item[(b)] There is $\beta>0$ and a unique $\eta\in A(a,\beta)$ such that
  $$
  x_u=x_{\sigma+t}=\widetilde{\varphi}_{\sigma+t}+\eta_t,\mbox{ where } u=\sigma+t, t\in[0,a].
  $$
  Moreover $\dot{x}(u)$ exists, if and only if, $\dot{\eta}(t)$ exists and $\dot{x}(u)=\dot{\eta}(t)$
\end{itemize}
\noindent
\textbf{Fact 3.2:}\\
The following are equivalent.
\begin{itemize}
  \item[(a)] $x$ is a solution of $(P)$.
  \item[(b)] $\exists\; a,\beta>0\;\exists\;\eta\in A(a,\beta):T\eta=\eta$
  (Indeed by Fact 3.1 we have for $u=\sigma+t\in(\sigma,\sigma+a)$
  $$
  \dot{x}(u)=f(u,x_u)=\dot{\eta}(t)=f(\sigma+t,\widetilde{\varphi}_{\sigma+t},\eta_t)\; ).
  $$
 \end{itemize}
Hence, the solutions of $(P)$ are in one to one correspondence with the fixed points of the operator $T$.
\\\\
\textbf{Fact 3.3:}\\
If the function $f(\tau,\psi)$, which defines $(P)$ is continuous with respect to $\psi$ for each fixed $\tau$ (In particular, if $f$ satisfies Caratheodory condition) and if $T(A(a,\beta))\subseteq A(a,\beta)$ then the operator $T$ has fixed points (That is, the problem $(P)$ has a solution on $[\sigma-r,\sigma+a]$). Also the above inclusion holds, in case that $f$ satisfies Caratheodory condition. (The proof of Fact 3.3 is contained in [H] pages 14,15).
\\\\
\textbf{Fact 3.4:}\\
The fixed points $\eta\in A(a,\beta)$ of the operator $T$ depend on the starting point $\sigma$ of solutions:
$$
\eta(t)=\eta(\sigma,t), t\in[0,a]
$$
This is evindend since although $\widetilde{\varphi}_{\sigma + s}+\eta_s=\widetilde{\varphi}_{\sigma'+s}+\eta_s$ for $\sigma\neq\sigma '$, where $\widetilde{\varphi}\in C_{\sigma,r}$ in the left hand side and $\widetilde{\varphi}\in C_{\sigma ',r}$ in the right hand side, generally $f(\sigma+s,\widetilde{\varphi}_{\sigma+s},\eta_s)\neq f(\sigma '+s,\widetilde{\varphi}_{\sigma '+s},\eta_s)$ and hence $T$ changes, if we change starting point for the solution.
\\\par
With the next proposition actually we see that under the hypothesis of continuous converge of a sequence $(f_k)_{k=1,2,...}$ to a function $f$, the requirement of uniform boundedness for $(f_k)_k$ in theorem 5.1 of [H] is redundant.
\\\\
\textbf{Proposition 3.5:}\\
Let $D$ be an open set in $\mathbb{R}\times C$ and $f_k: D\rightarrow \mathbb{R}$ be arbitrary functions, $k=0,1,2,...$. Suppose that
$$
f_k\xrightarrow{\mbox{cont}}f_0,\;\;\; k\rightarrow\infty
$$
(Hence $f_0$ is continuous by propositions 2.2. and 2.3).\\
If $x^{(0)}\in C_{\sigma,r}=C([\sigma-r,\sigma+a],\mathbb{R}^N)$ and $W=\{(t,x_t^{(0)}):t\in[\sigma,\sigma+a]\}\subseteq D$ then\\
$\exists\; V$ neighborhood of $W\;\exists\; M>0\;\exists\;k_0\in\mathbb{N}: |f_k(\tau,\psi)|\leq M$ for all $(\tau,\psi)\in V, k\geq k_0\;\& k=0$.
\\\\
\textbf{Proof:}
First we observe that $W$ is compact. Since by proposition 2.2 the sequence $(f_k)_{k=1,2,...}$ is exhaustive at each $(t,x_t^{(0)})\in W$ we get that
\begin{eqnarray}
\forall\; t\in[\sigma_0,\sigma_0+a_0]\;\exists\;\delta_t>0\;\exists\;k_t\in\mathbb{N}:\\
(\tau,\psi)\in S((t,x_t^{(0)});\delta_t)\Rightarrow|f_k(t,x_t^{(0)})-f_k(\tau,\psi)|<1, \mbox{ for } k\geq k_t .\nonumber
\end{eqnarray}
($S((t,x_t^{(0)});\delta_t)$ is the open ball with center $(t,x_t^{(0)})$ and radius $\delta_t$ in $\mathbb{R}\times C$).
If $(\tau,\psi)\in V$, then by (1) we get that
$$
|f_k(\tau,\psi)|\leq1+f_k(t_j,x_{t_j}^{(0)}),
$$
for some $j\in\{1,2,...,l\}$ and $k\geq k_t$.
Since $f_k(t_j,x_{t_j}^{(0)})\rightarrow f_0(t_j,x_{t_j}^{(0)}), k\rightarrow\infty$
for $j=1,2,...,l$ there are positive numbers $M_j$ such that
$$
|f_k(t_j,x_{t_j}^{(0)})|\leq M_j \mbox{ for } k\geq k_{tj}, j=1,2,..,l .
$$
We set
$$
M=\max_{\substack{1\leqq j\leqq l}}(M_j+1), k_0=\max_{\substack{1\leqq j\leqq l}} k_{t_j},
$$
hence it follows that
$$
|f_k(\tau,\psi)|\leqq M \mbox{ for } (\tau,\psi)\in V, k\geqq k_0.
$$
Finally by taking limit as $k\rightarrow\infty$ we get that
$$
|f_0(\tau,\psi)|\leqq M,  for (\tau,\psi)\in V.
$$
Before we state an improved version of the theorem on continuous depence of solutions of problem $(P)$ we make some necessary remarks.
\\\\
\textbf{Remarks 3.6:}
\begin{itemize}
  \item[(i)] In theorem 5.1 of [H] page 21, it is assumed that
  $$
  f_k(t,\psi)\rightarrow f_n(t,\varphi) \mbox{ as } k\rightarrow\infty, \psi\rightarrow\varphi
  $$
   that is $(f_k)_k$ converges continuously to $f_0$ only with respect the second variable for each dixed $t$. But in the proof (page 22) it is used continuous convergence with respect both variables. Namely,
  $$
  f_k(\sigma_k+s,\widetilde{\varphi}_{\sigma_k+s}^{(k)},\eta_s)\rightarrow f_0(\sigma_0+s,\widetilde{\varphi}_{\sigma_0+s}^{(0)},\eta_s), k\rightarrow\infty .
  $$
  If we wish to correct this inconsistency, we must take a number $\sigma >\sigma_0$ and define the operators $T_k, k=0,1,2,... $ as follows
  $$
  T_k\eta(t)=\left\{
                \begin{array}{ll}
                  \int_0^t f_k(\sigma+s,\widetilde{\varphi}_{\sigma+s}^{(k)},\eta_s)ds, t\in[0,\overline{a}](\overline{a}>0)\\
                  0, t\in[-r,0]
                \end{array}
              \right.
  $$
  Then the proof of the theorem leads to the fact that
  \begin{equation}
  \eta^{(k)}(\sigma,\cdot)\xrightarrow{||\;||}  \eta^{(0)}(\sigma,\cdot), k\rightarrow\infty ,
  \end{equation}
  where $\eta^{(k)}(\sigma,\cdot), k=0,1,...$ are the fixed points of the above operators $T_k$. On the other hand the conclusion of the theorem is the following.
  $$
  \eta^{(k)}(\sigma_k,\cdot)\xrightarrow{||\;||}\eta^{(0)}(\sigma_0,\cdot), \;k\rightarrow\infty
  $$
  (see also Fact 3.4), or equivalently,
  \begin{equation}
  \forall\;\varepsilon>0\exists k_1\in\mathbb{N}:|\eta^{(k)}(\sigma_k,t)-\eta^{(0)}(\sigma_0,t)|<\frac{\varepsilon}{3}\mbox{ for } k\geq k_1, t\in [0,\overline{a}].
  \end{equation}
  But by (2) we have that
  $$
  |\eta^{(k)}(\sigma,t)-\eta^{(0)}(\sigma,t)|<\frac{\varepsilon}{3}\mbox{ for } k \mbox{ large enough }, t\in [0,\overline{a}].
  $$
  Also, it is not hard to see that
  $$
  |\eta^{(0)}(\sigma,t)-\eta^{(0)}(\sigma_0,t)|<\frac{\varepsilon}{3}\mbox{ for } \sigma \mbox{ near } \sigma_0.
  $$
  Since,
  \begin{eqnarray*}
  |\eta^{(k)}(\sigma_k,t)-\eta^{(k)}(\sigma,t)|-|(\eta^{(k)}(\sigma,t)-\eta^{(0)}(\sigma,t))+(\eta^{(0)}(\sigma,t)-\eta^{(0)}(\sigma_0,t))|\leqq\\
  |\eta^{(k)}(\sigma_k,t)-\eta^{(0)}(\sigma_0,t)|<\frac{\varepsilon}{3},
  \end{eqnarray*}
we must have that
  $|\eta^{(k)}(\sigma_k,t)-\eta^{(k)}(\sigma,t)|<\varepsilon$ for $k$ large enough or $|\int_0^t f_k(\sigma_k+s,\widetilde{\varphi}_{\sigma_k+s}^{(k)},\eta_s)-f_k(\sigma+s,\widetilde{\varphi}_{\sigma+s}^{(k)},\eta_s)  ds|<\varepsilon,$ for $k$ large, $t\in [0,\overline{a}].$
  The above inequality is true in case that $|f_k(\sigma_k+s,\widetilde{\varphi}_{\sigma_k+s}^{(k)},\eta_s)-f_k(\sigma+s,\widetilde{\varphi}_{\sigma+s}^{(k)},\eta_s)|$ is arbitrarily small for $k$ sufficiently large. But this is the property of exhaustiveness for $(f_k)_k$ which in turn implies continuous convergence (See definition 2.1 and proposition 2.2). Hence continuous convergence of $(f_k)_k$ to $f_0$ with respect to both variables seems to be unavoidable.

  \item[(ii)] Another essential point, which must be restored is the following. In theorem 5.1 of [H] we consider arbitrary sequence $(x^{(k)})$ of solutions for the problem $(P_k):\dot{x}(t)=f_k(t,x(t))$, $x_{\sigma_k}=\varphi_k, k=0,1,2,...$. Also the problem $(P_k), k=1,2,...$, may have even infinite solutions.
  Suppose that each solution of $(P_k)$ has maximum domain the interval $(\sigma_k-r,\sigma_k+a_k)$, $a_k>0, k=1,2,...$. In the proof, it is shown that the corresponding operators $T_k$ of $(P_k)$ have fixed points on $A(\overline{a},\overline{\beta})$, where $\overline{a},\overline{\beta}$ are suitably choosen constants. The questions is that, how are we sure that the fixed points of $T_k$ fits with the considered solutions $x^{(k)}$? It may happened that $x^{(k)}$ is given by fixed points of $T_k$ on $A(\overline{a}_k,\overline{\beta}_k)$ with $\overline{a}_k< \overline{a}$. We restore this point with lemma 3.7
\end{itemize}

\textbf{Lemma 3.7:}
Suppose that the function $f$ defining problem $(P)$ satisfies the following $|f(\sigma+s,\widetilde{\varphi}_{\sigma+s},\eta_s)|\leq M$, for all $s\in[0,a],\eta\in A(a,\beta)$, where $a,\beta, M>0$, with $Ma<\beta$. Then, the operator $T$ corresponding to problem $(P)$ has fixed points on each $A(a_1,\beta)$ with $a_1\leq a$, which can be extended to fixed points of $T$ on $A(a,\beta)$.
\\\\
\textbf{Proof:}
The existence of fixed points for $T$ for each $a_1\leq a$ follows from our assumption on $f$ and fact 3.3.
Let $0<a_1<a$ and $\eta^{(1)}$ be a fixed point of $T$ on $A(a_1,\beta)$, that is
$$
\eta^{(1)}(t)=\int_0^t f(\sigma+s,\widetilde{\varphi}_{\sigma+s},\eta_s^{(1)}) ds, t\in[0,a_1]
$$
If,
$$
T_0\eta(t)=\int_0^t f((\sigma+a_1)+s,(\widetilde{\varphi_{\sigma+a_1}+\eta_{a_1}^{(1)}})+\eta_s)\; ds, t\in[0,a-a_1], \eta\in A(a-a_1,\beta)
$$
where $((\widetilde{\varphi_{\sigma+a_1}+\eta_{a_1}^{(1)}})\in C_{\sigma+a_1,r}$ and $\eta^{(0)}$ is a fixed point of $T_0$, we set
  $$
  \overline{\eta}(t)=\left\{
                \begin{array}{ll}
                  \eta^{(1)}(t), t\in [0,a_1]\\
                  (\widetilde{\eta^{(1)}_{a_1}})(t)+\eta^{(0)}(t-a_1),t\in [a_1,a], \mbox{ where } (\widetilde{\eta^{(1)}_{a_1}})\in C_{\sigma+a_1,r}
                \end{array}
              \right.
  $$
Then, for $t\in [0,a_1]$ we have that
$$
T\overline{\eta}(t)=\int_0^t f(\sigma+s,\widetilde{\varphi}_{\sigma+s},\overline{\eta}_s) ds=\int_0^t f(\sigma+s,\widetilde{\varphi}_{\sigma+s},\eta_s^{(1)}) ds=\eta^{(1)}(t)=\overline{\eta}(t)
$$
For $t\in[a_1,a]$ it holds the following
\begin{eqnarray*}
T\overline{\eta}(t)&=&\int_0^{a_1} f(\sigma+s,\widetilde{\varphi}_{\sigma+s},\overline{\eta}_s) ds+ \int^t_{a_1} f(\sigma+s,\widetilde{\varphi}_{\sigma+s},\overline{\eta}_s) ds\\
&=&\eta^{(1)}(a_1)+\int_0^{t-a_1} f(\sigma+a_1+u,\widetilde{\varphi}_{\sigma+a_1+u}+\overline{\eta}_{a_1+u}) ds\\
&=&(\widetilde{\eta^{(1)}(a_1)})(t)+\int_0^{t-a_1} f((\sigma+a_1)+u,(\widetilde{\widetilde{\varphi}_{\sigma+a_1}+\eta^{(1)}_{a_1}})_{\sigma+a_1+u}+\eta_u^{(0)}) du \\
& &\mbox{(By definition of }  \overline{\eta} \mbox{ we have that} \\
\overline{\eta}_{a_1+u}&=&(\widetilde{\eta_{a_1}^{(1)}})_{\sigma+a_1+u}+n^0_u\Rightarrow \widetilde{\varphi}_{\sigma+a_1+u}+\overline{\eta}_{a_1+u}=(\widetilde{\varphi_{\sigma+a_1}+\eta_{a_1}^{(1)}})_{\sigma+a_1+u}+\eta_u^0)\\\\
&=& (\eta_{a_1}^{(1)})(t)+\eta^{(0)}(t-a_1)\\\\
&=& \overline{\eta}(t))
\end{eqnarray*}
Hence
$$
T\overline{\eta}=\overline{\eta}
$$
and the proof of Lemma is completed.\\

\textbf{Theorem 3.8:}
Let $D$ be an open set in $\mathbb{R}\times C$, $f_k:D\rightarrow \mathbb{R}^N, k=0,1,2,...$ and $\{\varphi_k: k=0,1,2,...\}\subseteq C, \sigma_k\in\mathbb{R}, k=0,1,2,...$. We assume the following:
\begin{itemize}
  \item[(a)] $f_k$ satisfies Caratheodory condition on $D$ for $k=0,1,2,...$.
  \item[(b)] $f_k \xrightarrow{\mbox{cont}}f_0,k\rightarrow\infty$.
  \item[(c)] $\sigma_k\rightarrow\sigma_0, k\rightarrow\infty \& \varphi_k \xrightarrow{||\;||}\varphi_0, k\rightarrow\infty $
  \item[(d)] $x^{(k)}$ is a solution of the problem
  $$
    (P_k):\left\{
                \begin{array}{ll}
                  \dot{x}(t)=f_k(t,x(t))\\
                  x_{\sigma_k}=\varphi_k, k=0,1,2,...
                \end{array}
              \right.
  $$
\end{itemize}
and the solution $x^{(0)}$ of $(P_0)$ is unique. If $x^{(0)}$ has maximum domain the interval $[\sigma_0-r,\sigma_0+a_0), a_0>0$, then
\begin{itemize}
  \item $x^{(k)}$ have maximum domains the intervals $[\sigma_k-r,\sigma_k+a_k)$ with $a_k\geq a_0$, for $k$ sufficiently large
  \item $x^{(k)}\xrightarrow{\mbox{cont}}x^{(0)}, \; k\rightarrow\infty$\\
  (according with the generalized notion of continuous convergence, see definition 2.6). Hence $(x^{(k)})_k$ converges uniformly to $x^{(0)}$ on each common closed interval of their domains.
  \end{itemize}
\textbf{Proof:}
Let $W=\{(t,x_t^{(0)}):t\in[\sigma_0,\sigma_0+a^{'}]\}, a^{'}<a_0$. Then by proposition 3.5 there are open neighborhood $V$ of $W$, $M>0$ and $k_1\in\mathbb{N}$ such that
$$
|f_k(\tau,\psi)|\leqq M \mbox{ for all } (\tau,\psi)\in V, k\geqq k_1, k=0.
$$
We set
$$
\delta=d(W,V^c)>0
$$
$$
V_1=\{(\tau,\psi)\in D:d((\tau,\psi),W)<\frac{\delta}{2}\},
$$
where $d$ is the distance on $\mathbb{R}\times C$.
If $0<a<\frac{\delta}{2},0<\beta<\frac{\delta}{2}$, then
$$
(\tau,\psi)\in V_1, t\in[0,a],\eta\in A(a,b)\Rightarrow (\tau+t,\psi+\eta_t)\in V.
$$
According to our notations at the beginning of this paragraph, for each $\sigma\in[\sigma_0.\sigma_0 +a^{'}], x_\sigma^{(0)}\in C, (\widetilde{x_\sigma^{(0)}})\in C_{\sigma ,r},(\widetilde{x_\sigma^{(0)}})_{\sigma+t}\in C, \sigma_0\leq\sigma+t\leq\sigma_0+a^{'} $ and since $x_\sigma^{(0)}$ is uniformly continuous on $[\sigma_0, \sigma_0+a^{'}]$ we get that
\begin{equation}
\exists\;a^{''}>0, a^{''}<a, a^{'}:||(\widetilde{x_\sigma^{(0)}})_{\sigma+t}-x^{(0)}_{\sigma+t} ||<\frac{\beta}{3}\mbox{ for } \sigma\in[\sigma_0, \sigma_0+(a^{'}-a^{''})], t\in[0,a^{''}]
\end{equation}
(We note that $||(\widetilde{x_\sigma^{(0)}})_{\sigma+t}-x^{(0)}_{\sigma+t}||=\sup_{\substack{t\in[0,a^{''}]}}|x^{(0)}(\sigma)-x^{(0)}(\sigma+t)|$).
In particular for $\sigma=\sigma_0$ in (4) above taking into account that $\varphi^{(0)}=\widetilde{x^{(0)}_{\sigma_0}}\in C_{\sigma_0 ,r},\widetilde{\varphi}^{(k)}\in C_{\sigma_k ,r},\widetilde{\varphi}^{(k)}_{\sigma_k}=\varphi^{(k)}\xrightarrow{||\;||}\varphi^{(0)}=\widetilde{\varphi}^{(0)}_{\sigma_0},k\rightarrow\infty$ and $\sigma_k\rightarrow\sigma_0, k\rightarrow\infty$, we take the following
\begin{equation}
||\widetilde{\varphi}_{\sigma_0+t}^{(0)} -x^{(0)}_{\sigma+t} ||<\frac{\beta}{3}\mbox{ for } t\in[0,a^{''}].
\end{equation}
\begin{equation}
\exists\; k_2\in\mathbb{N}:||\widetilde{\varphi}_{\sigma_0+t}^{(0)} -\varphi^{(k)}_{\sigma_k+t} ||<\frac{\beta}{3}\mbox{ for } k\geq k_2, t>0.
\end{equation}
Now, since $||\widetilde{\varphi}_{\sigma_k+t}^{(k)}+\eta_t -x^{(0)}_{\sigma_0+t} ||\leqq ||\widetilde{\varphi}_{\sigma_k+t}^{(k)} -x^{(0)}_{\sigma_0+t} || + ||\eta_t||$, if we set $\overline{\beta}=\frac{\beta}{3}$, from all the above we get that
\begin{eqnarray*}
\exists\;\overline{a}>0\;\exists\; k_0\geqq k_1, k_2, (k_0\in\mathbb{N})&:&\\
\overline{a}\leqq a^{''}& &\; \& M \overline{a}<\overline{\beta}\; \& \\
(\sigma_k+t,\widetilde{\varphi}_{\sigma_k+t}^{(k)}+\eta_t)\in V,\mbox { for all } k\geqq k_0 &,& k=0, t\in[0,\overline{a}],\eta\in A(\overline{a},\overline{\beta})
\end{eqnarray*}
Hence,
\begin{eqnarray}
|T_k\eta(t)-T_k\eta(t^{'}) |&\leqq& \int_{t^{'}}^t|f_k(\sigma_k+s,\widetilde{\varphi}_{\sigma_k+t}^{(k)}+\eta_s)|ds\\\nonumber
&\leqq& M(t-t^{'}) \\ \nonumber
&\leqq& \overline{\beta}\mbox{, for all } 0\leqq t^{'}\leqq t\leqq\overline{a},k\leqq k_0, k=0,\nonumber
\end{eqnarray}
where $T_k$ is the corresponding operator to problem $(P_k)$
 $$
    \eta\in A(\overline{a},\overline{\beta})\rightarrow T_k\eta(t)=\left\{
                \begin{array}{ll}
                  \int_0^t f_k(\sigma_k+s,\widetilde{\varphi}_{\sigma_k+s}^{(k)}+\eta_s) ds, t\in [0,\overline{a}]\\
                  0, t\in[-r,0]
                \end{array}
              \right. ,k\in\mathbb{N}.
  $$
Since the solutions are considered with maximum demains, (7) and Lemma 3.7 implies that
\begin{eqnarray}
\mbox{the fixed points } \eta^{(k)} \mbox{ of the operators } T_k, \mbox{ corresponding to the solutions } x^{(k)}\\\nonumber
\mbox{ are defined on } [0,\overline{a}] \mbox{ for } k\geqq k_0, k=0.
\end{eqnarray}
Also from (7) and Arteza-Ascoli's theorem we get that
\begin{equation}
\{T_k\eta^{(k)}:k=1,2,...\} \mbox{ is relatively compact in } A(\overline{a},\overline{\beta})
\end{equation}
so,
\begin{equation}
\exists\;\eta^{\ast}A(\overline{a},\overline{\beta})\;\exists L\subseteq\mathbb{N}, L \mbox{ infinite: } \eta^{(k)}=T_k\eta^{(k)} \xrightarrow{||\;||}\eta^{\ast},k\rightarrow\infty, k\in L.
\end{equation}
Now, we set
$$
F_k: s\in[0,\overline{a}]\rightarrow f_k (\sigma_k+s,\widetilde{\varphi}_{\sigma_k+s}^{(k)},\eta_s^{(k)}), k=0,1,2,...
$$
Then by hypothesis $(b), (c)$ and (10) we get that
$$
F_k\xrightarrow{\mbox{cont}}F_0, k\rightarrow\infty, k\in L.
$$
Hence proposition 2.4 implies
$$
F_k\xrightarrow{||\;||}F_0, k\rightarrow\infty, k\in L.
$$
Let $\varepsilon >0$. We have
\begin{eqnarray*}
\sup_{\substack{t\in[0,\overline{a}]}}|\eta^{(k)}(t)-T_0^k(t)|& = & \sup_{\substack{t\in[0,\overline{a}]}} |\int_0^t F_k(s) ds-\int_0^t F_0(s) ds|\\
& \leqq & \sup_{\substack{t\in[0,\overline{a}]}}\int_0^t |F_k(s)-F_0(s)| ds   \\
& \leqq & \frac{\varepsilon}{\overline{a}}\cdot\overline{a}=\varepsilon,\mbox{ for } k \mbox{ sufficiently large, } k\in L.
\end{eqnarray*}
This means that
$$
\eta^{(k)}\xrightarrow{||\;||} T_0\eta^{\ast}
$$
and by (10) we have that
$$
T_0\eta^{\ast}=\eta^{\ast} .
$$
But the fixed point of $T_0$ is unique. So,
$$
\eta^{\ast}=\eta^{(0)}.
$$
Finally, since by (9) every subsequence of $\eta^{(k)}$ has a convergent subsequence, which according to the above reasoning converge to $\eta^{(0)}$ we take that
$$
\eta^{(k)}\xrightarrow{||\;||}\eta^{(0)}\rightarrow\infty, \mbox{ on } [0,\overline{a}]
$$
Since $\overline{a}$ depends only on the uniform continuity of $x^{(0)}$ on $[\sigma_0,\sigma_0+a^{'}]$, by successive stepping intervals of constant length $\overline{a}$ we take that
\begin{eqnarray*}
\eta^{(k)} \mbox{ are defined on } [0,a^{'}], \mbox{ for } k\geq k_0 \\
\eta^{(k)}\xrightarrow{||\;||}\eta^0\; \& \mbox{ on } [0,a^{'}], k\rightarrow\infty \mbox{ (see also [5]).}
\end{eqnarray*}
From the above we take the following
\begin{itemize}
  \item[(1st)] $a_k\geq a^{'}, k\geq k_0\Rightarrow a_k\geq a_0, k\geq k_0$
  \item[(2nd)] If $\theta_k\in[0,a_0), k=0,1,2,...$ and $\theta_k\rightarrow\theta ,k\rightarrow\infty$ then
  $$
  x^{(k)}(\sigma_k+\theta_k)=\widetilde{\varphi}^{(k)}(\sigma_k+\theta_k)+\eta^{(k)}(\theta_k)\rightarrow\widetilde{\varphi}^{(0)}(\sigma_0+\theta_0)+\eta^{(0)}(\theta_0)=x^{(0)}(\sigma_0+\theta_0), k\rightarrow\infty.
  $$
\end{itemize}
Hence $x^{(k)}\xrightarrow{\mbox{cont}}x^{(0)}, k\rightarrow 0$ and the proof os completed. $\Box$\\
A combination of the above theorem and a known result from Fourier series give us the following applications.
\\\\
\textbf{Applications 3.9:}
Suppose that $f:\mathbb{R}\rightarrow \mathbb{R}$ is a continuous periodic function of bounded variation, with period, say for simplicity, $2\pi$. Then it is known that the sequence of Fourier partial sums of $f$
$$
S_n(x)=\frac{a_0}{2}+\sum_{k=1}^n(a_k\cos kx+b_k\sin kx), n\in\mathbb{N},
$$
converges continuously to $f$
$$
(a_k=\frac{1}{\pi}\int_0^{2\pi} f(t)\cos kt dt,\;\; b_k=\frac{1}{\pi}\int_0^{2\pi} f(t)\sin kt dt, k=0,1,2,..., \mbox{ see [Z]})
$$
Hence, we can approximate uniformly the solution of the problem
$$
 (P) \left\{
                \begin{array}{ll}
                  \dot{x}(t)=f(x(t))\\
                  x(0)=c.
                \end{array}
              \right.
$$
by the solutions of the problems
$$
 (P_n) \left\{
                \begin{array}{ll}
                  \dot{x}(t)=S_n(x(t))\\
                  x(0)=c.
                \end{array}
              \right.,\;\; n=1,2,...
$$
We close with the following question:\\
can we prove theorem 3.8 under the assumption of continuous convergence of the sequence $(f_k)$ to function $f_0$ only with respect to the second variable (adding some more assumptions, if necessary), or we can found a counterexample which verifies that the theorem is false in this case?

\end{document}